\documentclass[11pt, oneside]{article}
\usepackage{amsmath,amsfonts,amssymb,amscd,theorem}
\date{}

\newtheorem{thm}{Theorem}
\newtheorem{prop}{Proposition}
\newtheorem{lem}{Lemma}

\begin{document}

\title{Harmonic metrics on  unipotent bundles over  quasi-compact K\"ahler manifolds}

\author{J\"urgen Jost\footnote{Max Planck Institute for Mathematics in the Sciences,
Leipzig, Germany;}, Yi-Hu Yang\footnote{Department of
Mathematics, Tongji University, Shanghai, China;} \thanks{Supported partially
by NSF of China (No.10771160);}, and Kang Zuo\footnote{Department of Mathematics,
Mainz University, Mainz, Germany;} \thanks{Supported partially by SFB/TR 45 Periods,
Moduli Spaces and Arithmetic
of Algebraic Varieties of the DFG.}}

\maketitle

\abstract{In this note, we propose an approach to the study of the analogue 
for unipotent harmonic bundles of Schmid's Nilpotent Orbit Theorem.
Using this approach, we construct harmonic metrics on unipotent bundles over
quasi-compact K\"ahler manifolds with carefully controlled
asymptotics near the compactifying divisor; such a metric is unique up to some isometry.
Such an asymptotic behavior is canonical in some sense.}

\section{Introduction}
In this paper, we construct a harmonic metric  on a unipotent bundle  over a
quasi-compact K\"ahler manifold. A quasi-compact K\"ahler manifold
here is a K\"ahler manifold that can be compactified to a complex
variety by adding a divisor that consists of smooth hypersurfaces with
at most normal crossings as singularities. The key point of our
construction is the prescribed and
carefully controlled asymptotic behavior of the harmonic metric when
approaching that compactifying divisor. This control is canonical in some sense.
In fact, we have shown existence results for harmonic metrics on unipotent bundles
already in our earlier papers \cite{jost-zuo}, but the
methods employed there do not yield the control near the compactifying
divisor. Therefore, here we shall develop a new and more subtle method.

We now introduce some notation. Let $\overline X$ be an $n$-dimensional compact
K\"ahler manifold, $D$ a normal crossing divisor; $X=\overline X\setminus D$ then is
the quasi-compact K\"ahler manifold that we are interested in.
Let $\rho: \pi_1(X)\to SL(r, \mathbb{Z})\subset SL(r, \mathbb{C})$ be a linear
representation. Equivalently, we have a flat vector bundle $L_\rho$ over $X$.

Take $p\in D$ and a small neighborhood in $\overline X$ of $p$, say
$\Delta^n$, where $\Delta$ is the unit disk; by our normal crossing assumption,
we can assume that $\Delta^n\cap X=(\Delta^*)^k\times\Delta^{n-k}$,
where $\Delta^*$ is the punctured disk.
Let $(z_1, z_2, \cdots, z_n)$ be the corresponding local complex
coordinate at $p$ covering $\Delta^n$, $z_i=r_ie^{{\sqrt{-1}}\theta_i}
(1\le i\le n)$. Let $\gamma_1, \gamma_2, \cdots, \gamma_k$ be the
generators of $\pi_1((\Delta^*)^k\times\Delta^{n-k})$, where
$\gamma_i$ corresponds to the $i$-th $\Delta^* (1\le i\le
k)$. We restrict $\rho$ to $\pi_1((\Delta^*)^k\times\Delta^{n-k})$
(note that the induced map from
$\pi_1((\Delta^*)^k\times\Delta^{n-k})$ into $\pi_1(X)$ is not
necessarily injective), and denote $\rho(\gamma_i)$ still by
$\gamma_i$. We note that $\gamma_1, \gamma_2, \cdots, \gamma_k$
commutate, which will be important for our constructions below.

In this note, we always assume that {\it each $\gamma_i$ is
unipotent}; so $N_i=\log\gamma_i$ is nilpotent and the
$N_1, \cdots, N_k$  commutate. We call $L_\rho$ a {\it unipotent
bundle}. If $L_\rho$ is endowed with a harmonic metric $h$ with {\it
tame} growth condition (Simpson's terminology \cite{simpson}) near the divisor,
we call $(L_\rho, h)$ a {\it unipotent harmonic bundle}.
For the notion of general harmonic bundles and its history and related definitions,
we refer the reader to \cite{simpson1}.

The notion of a unipotent harmonic bundle is a natural generalization
of a variation of Hodge structures introduced by P. Griffiths \cite{griffiths, simpson}.
In the study of variations of Hodge structures,
it is very important to understand the asymptotic
behavior near the divisor of the Hodge metric \cite{cattani-kaplan-schmid2, kashiwara-kawai};
likewise, for the study of a unipotent harmonic bundle, we also need to have
a good understanding of the asymptotic behavior of the endowed tame harmonic metric.
In this aspect, Mochizuki \cite{mochizuki}, along the way initiated by
Cattani-Kaplan-Schmid \cite{cattani-kaplan-schmid}, provides
a systematic description for the asymptotic behavior of tame harmonic
metrics on a unipotent bundle.
In this paper, we try to
propose a slightly different approach to the study of asymptotic behavior of a unipotent harmonic bundle;
our approach is perhaps able to be considered as
an analogue of Schmid's Nilpotent Orbit Theorem \cite{schmid}.

In order to formulate the idea, we first go back to Schmid's earlier paper \cite{schmid};
there, in the case of one variable, an equivalent description of the asymptotic behavior of the
period mapping (equivalently, the Hodge metric) is developed in terms
of an equivariant geodesic embedding of the upper half plane into the
period domain---a homogeneous (complex) manifold; although this was
not explicitly stated, it clearly is implied by the nilpotent orbit theorem
and $SL_2$-orbit theorem for one variable. The existence and
uniqueness of such a geodesic embedding are obtained by using
Jacobson-Morosov's theorem and a result of Kostant \cite{kostant}. We
hope to generalize this description to the case of several variables;
consequently, we also need to extend Jacobson-Morosov's theorem and
the result of Kostant in some sense. However, we will not directly use
such a generalization to obtain general information about the
asymptotic behavior of a unipotent harmonic bundle; instead we first use
this idea to construct an initial metric of finite energy on a
unipotent bundle and then deform this initial metric to a harmonic
metric of finite energy without changing the behavior near the divisor.
So, our tame harmonic metric is a more special one, namely one with trivial filtered
(or parabolic) structure in the sense of local systems \cite{simpson, mochizuki};
a general tame harmonic metric may be of infinite energy \cite{jyz}.

The idea of using a geodesic embedding of the upper half plane as the
asymptotic behavior of a unipotent harmonic bundle in the case of one
variable has already been successfully used in our previous work
\cite{jost-yang-zuo2}---the study of cohomologies for unipotent harmonic bundles
over a noncompact curve, more precisely, the study of the $L^2$-Poincar\'e lemma
and the $L^2$-Dolbeault lemma for unipotent harmonic bundles.

In order to apply Schmid's description of a geodesic embedding of the upper half plane
as the asymptotic behavior of a Hodge metric to harmonic metrics in the higher
dimensional case, one needs to suitably extend Jacobson-Morosov's theorem and the
result of Kostant; we feel that such an extension is an appropriate (algebraic)
substitute in the harmonic bundles theory of Cattani-Kaplan-Schmid's $SL_2$-orbit
theorem for several variables in the Hodge theory.
In \S 2, we work out this extension; mainly, we use the notions of
parabolic subgroups (algebras) and the corresponding horospherical
decomposition. Our argument is slightly geometric; we think that there
should also exist a purely Lie-theoretic proof.

In Schmid's description, one of the key points is how to get a related geodesic
embedding by using one single nilpotent element; this is achieved by finding a
semi-simple element by using Jacobson-Morosov's theorem, then Kostant's result
implies uniqueness in a certain sense.
For several commuting nilpotent elements $N_1, \cdots, N_k$ in the
present setting, correspondingly we hope to get a group of semi-simple
elements, $Y_1, \cdots, Y_k$, {\it which are commutative} and such
that each pair $(N_i, Y_i)$ can be extended to an $\mathfrak{sl}_2$-triple.
(We here remark that in general these triples are not commutative.)
The commutativity of $Y_1, \cdots, Y_k$ implies that they are contained in
a maximal abelian subspace; this motivates us to use the theory of parabolic
subgroups (subalgebras). We also remark that these semisimple elements are
unique after fixing a maximal abelain subspace.

After getting the semi-simple elements $Y_1, \cdots, Y_k$, in \S 3, we
are able to construct initial metrics of finite energy. These metrics
are not yet harmonic, but for our local construction, it is indeed
harmonic when restricted to a punctured disk transversal to the divisor; this is enough
for us to prove the asymptotic behavior of the obtained harmonic metric.

Let $\mathcal{P}_r$ be the set of positive
definite hermitian symmetric matrices of order $r$ with determinant $1$. $SL(r,
\mathbb{C})$ acts transitively on $\mathcal{P}_r$ by
\[
g\circ H=:{g}H{\bar g}^t, ~H\in\mathcal{P}_r, g\in SL(r,
\mathbb{C}).
\]
Obviously, the action has the isotropy subgroup $SU(r)$ at the
identity $I_r$. Thus $\mathcal{P}_r$ can be identified with the
symmetric space of noncompact type $SL(r, \mathbb{C})/SU(r)$, and can be uniquely endowed with
an invariant metric up to some constant.
Our local construction at $p$ then takes the following form
\begin{equation}
\label{H0}
H_0(z_1, z_2, \cdots, z_n)=\exp({\frac 1{2\pi}}\sum_{i=1}^k\theta_iN_i)
\circ\exp(\sum_{i=1}^k({\frac 1 2}\log|\log r_i|)Y_i).
\end{equation}
Such a construction is compatible for all $z_i$-directions; namely, for each
$z_i$-direction, it gives an equivariant geodesic embedding of the upper half plane.

In \S 4, we then deform the initial metric of finite energy to a harmonic metric
with the same asymptotic behavior. In order to make this deformation successful,
one needs to impose a geometric condition on the representation $\rho$, namely
semi-simplicity; such a condition first appeared in \cite{donaldson} and then
\cite{corlette}. Under this condition, we can deform the initial metric to a harmonic
one of finite energy. The harmonic metric is pluriharmonic by using Siu's Bochner
technique for harmonic map theory into a Riemannian manifold
due to Sampson \cite{sampson}. The essential difficulty is to prove that the harmonic
metric has the same asymptotic behavior as the initial metric; to this end,
restricting the harmonic metric and the above local construction of the initial
metric to a punctured disk transversal to the divisor so that both of them are
harmonic, we can then show that the distance function between both metrics is
actually bounded on the punctured disk, which accordingly implies that the harmonic
metric and the initial metric have the same asymptotic behavior.

We can now state our main results.

\begin{thm} Let $\overline X$ be an $n$-dimensioanl
compact K\"ahler manifold, $D$ a normal crossing divisor; set
$X=\overline X\setminus D$. Let $\rho: \pi_1(X)\to SL(r,
\mathbb{Z})\subset SL(r, \mathbb{C})$ be a linear semisimple
representation that is unipotent near the divisor, and $L_\rho$ the
corresponding unipotent bundle. Then, up to a certain isometry in $\mathcal{P}_r$,
there exists uniquely a harmonic metric
on $L_\rho$ with the same asymptotic behavior as $H_0$ (see {\rm (\ref{H0})})near the divisor.
\end{thm}

Our proof in \S4 also gives the following corollary where we do not assume that $\rho$
is semisimple, which is perhaps able to be considered as an analogue of Schmid's Nilpotent
Orbit Theorem. Based on this, we consider the asymptotic behavior {\rm (\ref{H0})}
as a canonical one.

\begin{thm} Any harmonic metric $H$ on $L_\rho$ of finite energy must have the same
asymptotic behavior as the local construction {\rm (\ref{H0})} near the divisor; namely
the distance function ${\text{dist}}_{\mathcal{P}_r}(H, H_0)$ is bounded near the divisor.
\end{thm}

\noindent {\bf Remark:} {\it The point in Theorem 2 that is unsatisfactory is that we can not yet give
the decay estimate of the distance function ${\text{dist}}_{\mathcal{P}_r}(H, H_0)$
near the divisor presently.}

\vskip .3cm
The research of harmonic metrics on noncompact manifolds was initiated by
Simpson \cite{simpson} in the complex one dimensional case under a more
algebraic geometric background; there, he also suggest that one should study
the analogues for harmonic bundles of the nilpotent orbit theorem and
the $SL_2$-orbit theorem in the Hodge theory; this paper can be considered as an attempt to it.
The general construction of harmonic metrics on a quasi-compact
K\"ahler manifold was later considered in \cite{jost-zuo, jost-zuo2} in a more
general setting --- equivariant harmonic maps.
\cite{jost-zuo} also obtained a harmonic metric of finite
energy; however, due to their construction, the
behavior at infinity of the metric is not controlled; in particular,
there is no norm estimate for flat sections when translated into the case of harmonic bundles.
Hopefully, our construction will give a new understanding of asymptotic behavior
of variation of Hodge structures or period mapping \cite{cattani-kaplan-schmid}.

This work was begun in February-April of 2004 when the second author was visiting CUHK and HKU in Hong Kong and continued in April-June of 2007 when he was visiting the Universities of Mainz and Essen; its main part
was finished when he visited the Max-Planck Institute for Mathematics in the Sciences, Leipzig in January-February and July-August of 2009. He thanks the above institutes for hospitality and good working conditions; especially, he wants to thank Professors Ngming Mok, H\'el\`ene Esnault and Eckart Viehweg for their kind invitation.

\section{Construction of semisimple elements}
In this section, we first recall Morosov-Jacobson's theorem, the theorem of Kostant and the horospherical decomposition associated to a parabolic subgroup; then we present the construction of semisimple elements $Y_1, \cdots, Y_k$.
\subsection{Some Lie-theoretic and geometric preliminaries}
{\bf 2.1.1. Morosov-Jacobson's theorem}
(cf. e.g. \cite{kostant}): {\it Let $G$ be a noncompact real simple Lie group, $\mathfrak{g}$ the corresponding Lie algebra. Assume $N$ is a nilpotent element in $\mathfrak{g}$. Then, one can extend $N$ to an $\mathfrak{sl}(2, \mathbb{R})$-embedding into $\mathfrak{g}$: $\{N, Y, N^-\}\subset\mathfrak{g}$ satisfying
\[
[N, Y]=2N, ~~[N, Y^-]=-2N^-, ~~[N, N^-]=Y.
\]}

Geometrically, this means the following: Let $X$ be the corresponding
Riemannian symmetric space of $G$. Then through any fixed point, there
exists a geodesic embedding into $X$ of the upper half plane whose horocycles are generated by $N$ and whose geodesics perpendicular to the horocycles are generated by the corresponding semi-simple element (or say the orbits of the corresponding one-parameter group).
\\
\\
{\bf 2.1.2. Kostant's theorem} (cf. \cite{kostant}): {\it Let $\mathfrak{g}$ and $N$ be as in Morosov-Jacobson's theorem. Then,
\\
1) the elements in $Im({\text{ad}}N)\cap Ker({\text{ad}}N)$ are nilpotent;
\\
2) if $\{N, Y, N^-\}$ is an $\mathfrak{sl}_2$-embedding extended by $N$ in $\mathfrak{g}$ as in Morosov-Jacobson's theorem, satisfying
\[
[N, Y]=2N, ~~[N, Y^-]=-2N^-, ~~[N, N^-]=Y,
\]
then $Y$ is unique up to nilpotent elements in $Im({\text{ad}}N)\cap Ker({\text{ad}}N)$.}

\vskip .3cm
Consequently, we have: Let $\mathfrak{g}=\mathfrak{k}+\mathfrak{p}$ be the Cartan decomposition and $\mathfrak{a}$ a maximal abelian subalgebra in $\mathfrak{p}$. If $Y$ lies in $\mathfrak{a}$, then, such a $Y$ is unique. Namely, for a fixed maximal abelian subalgebra $\mathfrak{a}\subset\mathfrak{p}$, $Y\in\mathfrak{a}$ is unique.

Geometrically, this means that, through a fixed point in $X$, the geodesic embedding in X of the upper half plane corresponding to $N$ in Morosov-Jacobson's theorem is unique.
\\
\\
{\bf 2.1.3. Parabolic subgroups (subalgebras) and the corresponding horospherical decompositions} (cf. e.g. Borel and Ji's book \cite{borel-ji})
\\
Let $\mathfrak{g}$ be a real semi-simple Lie algebra of noncompact type, $G$ the corresponding noncompact semisimple Lie group, $\mathfrak{g}=\mathfrak{k}+\mathfrak{p}$ Cartan decomposition, and $\mathfrak{a}$ a maximal abelian subalgebra of $\mathfrak{p}$. Relative to $\mathfrak{a}$, one has the restricted root system, denoted by $\Phi$; choose a simple root system, denoted by $\Delta$; denote the corresponding positive root system by $\Phi^+$.

Let $I$ be a subset of $\Delta$. One can construct a corresponding parabolic subalgebra $\mathfrak{p}_I$ (resp. $P_I$) of $\mathfrak{g}$ (resp. $G$), called the standard parabolic subalgebra (subgroup) relative to $I$. Any parabolic subalgebra (subgroup) is conjugate to such a one under $G$ and also under $K$, the maximal compact subgroup of $G$ corresponding to $\mathfrak{k}$; moreover, for any two distinct subsets $I, I'$ of $\Delta$, $P_I, P_{I'}$ are not conjugate under $G$.

Let $\mathfrak{m}$ be the centralizer in $\mathfrak{k}$ of $\mathfrak{a}$, $\mathfrak{a}_I=\cap_{\alpha\in I}Ker\alpha\subset\mathfrak{a}$, $\mathfrak{a}^I$ the orthogonal complement in $\mathfrak{a}$ of $\mathfrak{a}_I$; let $\Phi^I$ be the set of roots generated by $I$.
Set
\[
\mathfrak{n}_I=\sum_{\alpha\in\Phi^+-\Phi^I}\mathfrak{g}^\alpha,~~~
\mathfrak{m}_I=\mathfrak{m}\oplus\mathfrak{a}^I\oplus\sum_{\alpha\in\Phi^I}\mathfrak{g}^\alpha.
\]
Then, $\mathfrak{p}_I=\mathfrak{n}_I\oplus\mathfrak{a}_I\oplus\mathfrak{m}_I$ is the desired parabolic subalgebra, its corresponding subgroup of $G$ denoted by $P_I$.

Corresponding to the above decomposition of $\mathfrak{p}_I$, one also has the decomposition of $P_I$ in the level of groups: Let $N_I, A_I, M_I$ be the Lie subgroups of $G$ having the Lie algebras $\mathfrak{n}_I, \mathfrak{a}_I, \mathfrak{m}_I$ respectively, then
\[
P_I=N_IA_IM_I\cong N_I\times A_I\times M_I,
\]
where by $\cong$ we mean an analytic diffeomorphism, i.e. the map
$$
(n, a, m)\to nam\in P_I ~~(n\in N_I, a\in A_I, m\in M_I)
$$
is an analytic diffeomorphism. This is the so-called {\it Langlands decomposition} of $P_I$.
When $I$ is empty, $P_I=P_\emptyset$ is a minimal parabolic subgroup of $G$, the corresponding decomposition is $P_\emptyset=N_\emptyset AM$, here $A=\exp\mathfrak{a}, M=\exp\mathfrak{m}$.

Let $K$ be the maximal compact subgroup with the Lie algebra $\mathfrak{k}$ of $G$. The Iwasawa decomposition $G=N_\emptyset AK$ tells us that $G=P_IK$ for any subset $I$ of $\Delta$. So, $P_I$ acts transitively on the symmetric space  $X=G/K$. Thus, the Langlands decomposition for $P_I$ induces a decomposition of $X$ associated to $P_I$, called the {\it horospherical decomposition}
\[
X\cong N_I\times A_I\times X_I,
\]
where $X_I=M_I/(M_I\cap K)$, called the boundary symmetric space associated to $P_I$; by $\cong$ we again mean an analytic diffeomorphism, i.e. the map
\[
(n, a, m(M_I\cap K))\to namK\in X ~~(n\in N_I, a\in A_I, m\in M_I)
\]
is an analytic diffeomorphism.

Let $P$ be a parabolic subgroup of $G$. As mentioned before, it is conjugate to a unique standard parabolic subgroup $P_I$ under $K$. Choose $k\in K$ such that under $k$, $P$ is conjugate to $P_I$, denoted by $P=^kP_I$. Define
\[
N_P=^kN_I, ~~A_P=^kN_I, ~~M_P=^kM_I.
\]
Though the choice of $k$ is not unique, the subgroups $N_P, A_P, M_P$ are well-defined. We call $N_P, A_P$ (resp. the corresponding Lie algebras $\mathfrak{n}_P, \mathfrak{a}_P$) the {\it unipotent radical, the split component} of $P$ (resp. $\mathfrak{p}$) respectively.
Thus, we can translate the Langlands decomposition of $P_I$ into that of $P$, namely
\[
P=N_PA_PM_P\cong N_P\times A_P\times M_P.
\]
Consequently, we also have the horospherical decomposition of $X=G/K$ associated to $P$
\[
X\cong N_P\times A_P\times X_P,
\]
where $X_P=M_P/(M_P\cap K)$, called the boundary symmetric space associated with $P$.

\subsection{The construction of semisimple elements}
Now, we return to the setting of \S 1. As mentioned there, $N_1, N, \cdots, N_k$ are some commutative nilpotent matrices. By Engel's theorem, we can assume that all of them are upper triangular with the entries of the diagonal being zero, and hence that $\gamma_1, \gamma_2, \cdots, \gamma_k$ are upper triangular with the entries of the diagonal being $1$.

Choose a maximal parabolic subalgebra (resp. subgroup) $\mathfrak{p}$ (resp. $P$) of $\mathfrak{sl}(r, \mathbb{R})$ (resp. $SL(r, \mathbb{R})$) the unipotent radical of which contains $N_1, N, \cdots, N_k$ (resp. $\gamma_1, \gamma_2, \cdots, \gamma_k$); furthermore, we can choose $\mathfrak{p}$ (resp. $P$) so that its split component is contained in the set of diagonal matrices of $\mathfrak{sl}(r, \mathbb{R})$ (resp. $SL(r, \mathbb{R})$).

We remark that the set of diagonal matrices in $\mathfrak{sl}(r, \mathbb{R})$ is a maximal abelian subspace contained in the noncompact part of a Cartan decomposition of $\mathfrak{sl}(r, \mathbb{R})$; and that the use of the Engel's theorem shows that we consider the set of diagonal matrices in $\mathfrak{sl}(r, \mathbb{R})$ as such a maximal abelian subspace.
We also remark that such a choice of parabolic subalgebras (resp. subgroups) is not unique, even for a fixed maximal abelian subspace.

Let $P=N_PA_PM_P$ ($\mathfrak{p}=\mathfrak{n}_P\oplus\mathfrak{a}_P\oplus\mathfrak{m}_P$, $N_P=\exp\mathfrak{n}_P, A_P=\exp\mathfrak{a}_P$) be the Langlands decomposition of $P$, correspondingly, we have the horospherical decomposition
\[
SL(r, \mathbb{R})/SO(r):=X=N_PA_PX_P,
\]
where $X_P=M_P/(M_P\cap SO(r))$, the boundary symmetric space associated with $P$.
From the previous choice for $\mathfrak{p}$ (resp. $P$), we know that $\mathfrak{a}_P$ (resp. $A_P$), as a set of matrices, is contained in the set of diagonal matrices.

For $N_i\in\mathfrak{n}_P$, $i=1, \cdots, k$, by the Morosov-Jacobson's theorem, we can extend it to an embedding into $\mathfrak{sl}(r, \mathbb{R})$ of $\mathfrak{sl}(2, \mathbb{R})$, say $\{N_i, Y_i, N_i^-\}\subset\mathfrak{sl}(r, \mathbb{R})$, satisfying
\[
[N_i, Y_i]=2N_i, ~~[N_i, Y_i^-]=-2N_i^-, ~~[N_i, N_i^-]=Y_i.
\]

\begin{lem}
For this embedding, one can choose $Y_i$ such that it lies in $\mathfrak{a}_P$, and hence such a semisimple element is unique.
\end{lem}

In the following, we will show that this can actually be done by using (the geometric interpretations of) the Morosov-Jacobson theorem and the Kostant theorem, and the horospherical decomposition associated to $P$.
It should be an interesting question whether one can give  a purely
Lie-theoretic proof of this result.

Using the Killing form, we can easily show that the factors
$\mathfrak{n}_P, \mathfrak{a}_P, \mathfrak{m}_P$ of $\mathfrak{p}$ are
orthogonal to each other, so the orbits of $N_P$ and $A_P$ in the horospherical decomposition are orthogonal and also orthogonal to the boundary symmetric space $X_P$; in addition, we also can consider $X_P$ as the set of fixed points at infinity of $N_P$ and $A_P$.

Fix a point $x_0\in X$. By $\infty_0$, we denote the (unique) intersection point of the orbit of $x_0$ under $\exp(tN_i), t\in\mathbb{R}$ (denoted by $\exp(tN_i)\circ x_0$) with $X_P$. One has then a unique geodesic in $X$ connecting $x_0$ and $\infty_0$, denoted by $\sigma_0$. Clearly, the orbit $\exp(tN_i)\circ\sigma_0, t\in\mathbb{R}$ is a geodesic embedding into $X$.
We need to show that the orbit $\exp(tN_i)\circ\sigma_0, t\in\mathbb{R}$ is a geodesic embedding of the upper half plane and hence the orbit of any point in it under $\exp(tN_i)$ is its horocycle.
This is a consequence of the Morosov-Jacobson theorem and the Kostant theorem. By
the Morosov-Jacobson theorem, through $x_0$, we have a geodesic embedding of the upper half plane into $X$; the Kostant theorem implies that such an embedding is unique. On the other hand, $\exp(tN_i)\circ x_0$ is contained in this embedding, so the intersection point of this embedding with $X_P$ is also $\infty_0$. That is to say, this embedding also contains the geodesic $\sigma_0$. So, this embedding is just $\exp(tN_i)\circ\sigma_0, t\in\mathbb{R}$.
Furthermore, by the horospherical decomposition associated to $P$, we
can choose a semisimple element $Y_i$ in $\mathfrak{a}_P$, the orbits
in this embedding of the one-parameter group of which are geodesics
perpendicular to the horocycles. Such a $Y_i$ is just the desired
one. This completes the proof of the Lemma.
\\
\\
{\bf Remarks} 1) {\it From the construction, $Y_1, \cdots, Y_k\in\mathfrak{a}_P$ seem to depend on the choice of the parabolic subalgebra $\mathfrak{p}$. But, by  Kostant's theorem, we know that $Y_1, \cdots, Y_k$ only depend on the choice of a maximal abelian subspace in the noncompact part of a Cartan decomposition; namely they are unique up to some conjugations when we require that all of them be contained in a maximal abelian subspace.
So, when we fix a maximal abelian subspace $\mathfrak{a}$, we can uniquely get semisimple elements $Y_1, \cdots, Y_k\in\mathfrak{a}$ such that $[N_i, Y_i]=2N_i, i=1, \cdots, k$.}
2) {\it The above construction for semisimple elements works for any semisimple Lie algebra, not only for $\mathfrak{sl}(r, \mathbb{R})$; so, this may provide a new way to understand variations of Hodge structures and their degeneration, different from Cattani-Kaplan-Schmid's theory \cite{cattani-kaplan-schmid}.}

\section{The construction of initial metrics (maps) and their asymptotic behavior}
As in the Introduction, let $\overline X$ be a compact K\"ahler manifold, $D$ a normal crossing divisor, $X=\overline X\setminus D$. Taking $p\in D$ and a small neighborhood $U$ at $p$, then $X\cap U$ is of the form $(\Delta^*)^k\times\Delta^{n-k}$. Let $(z_1, z_2, \cdots, z_n)$ be a local complex coordinate on $U$ with
\[
(\Delta^*)^k\times\Delta^{n-k}=\{(z_1, z_2, \cdots, z_n): z_1\neq 0, z_2\neq 0, \cdots, z_k\neq 0\}.
\]
On $(\Delta^*)^k\times\Delta^{n-k}$, one has the following product metric
\[
ds_P^2={\frac{\sqrt{-1}}2}\big[\sum_{i=1}^k{\frac{dz_i\wedge d{\overline z}_i}{|z_i|^2(\log|z_i|)^2}}+\sum_{i=k+1}^ndz_i\wedge d\overline z_i\big].
\]

In general, one has the following
\begin{prop}
There exists a complete, finite volume K\"ahler metric on $X$ which is quasi-isometric to the metric of the above form near any point in the divisor $D$.
\end{prop}
\emph{Proof.} cf. \cite{cornalba-griffiths}.

\vskip .3cm
Due to the above proposition, when we consider local constructions and
various estimates near the divisor in this section, if not involving
the derivatives of the K\"ahler metric (in fact we indeed do not need
estimates involving  the derivatives), we always use the above local product metric $ds^2_P$.

\subsection{The construction of initial metrics of finite energy}
Using the previous constructions for semisimple elements $Y_1, \cdots, Y_n$, we can construct an initial map from the universal covering of $(\Delta^*)^k\times\Delta^{n-k}$ into $SL(r, \mathbb{C})/U(r)$, which is $\rho$-equivariant; equivalently, we can also consider such a map as a metric on the corresponding flat bundle ${L_\rho}|_{(\Delta^*)^k\times\Delta^{n-k}}$.

To this end, let us first give some
preliminaries. Let $\mathcal{P}_r$ be the set of positive
definite hermitian symmetric matrices of order $r$ with determinant $1$. $SL(r,
\mathbb{C})$ acts transitively on $\mathcal{P}_r$ by
\[
g\circ H=:{g}H{\bar g}^t, ~H\in\mathcal{P}_r, g\in SL(r,
\mathbb{C}).
\]
Obviously, the action has the isotropy subgroup $SU(r)$ at the
identity $I_r$. Thus $\mathcal{P}_r$ can be identified with the
coset space $SL(r, \mathbb{C})/SU(r)$, and can be uniquely endowed with
an invariant metric{\footnote{In terms of matrices, such an
invariant metric can be defined as follows. At the identity $I_r$,
the tangent elements just are hermitian symmetric matrices of trace-free; let $A, B$ be
such matrices, then the Riemannian inner product $<A,
B>_{\mathcal{P}_r}$ is defined by $tr(AB)$. In general, let $H\in
\mathcal{P}_r$, $A, B$ two tangent elements at $H$, then the
Riemannian inner product $<A, B>_{\mathcal{P}_r}$ is defined by
$tr(H^{-1}AH^{-1}B)$.}} up to some constants. In particular, under
such a metric, the geodesics through the identity $I_r$ are of the
form $\exp(tA)$, $t\in \mathbb{R}$, $A$ being a hermitian symmetric matrix of trace-free.
\\
\\
{\bf 3.1.1 Local construction}
Let $(z_1, z_2, \cdots, z_n)$ be the above complex coordinates on $(\Delta^*)^k\times\Delta^{n-k}$, $z_i=r_ie^{{\sqrt{-1}}\theta_i}$, $0<r_i<1, -\infty<\theta_i<\infty$.
Set
\[
H_0(z_1, z_2, \cdots, z_n)=\exp({\frac 1{2\pi}}\sum_{i=1}^k\theta_iN_i)\circ\exp(\sum_{i=1}^k({\frac 1 2}\log|\log r_i|)Y_i),
\]
which is independent of $z_{k+1}, \cdots, z_n$.
Clearly, it is $\rho$-equivariant.
Similar to \cite{jost-yang-zuo}, one can show that under the product metric $ds^2_P$,
$H_0$ has finite energy
{\footnote{We here remark that, in Proposition 1 of \cite{jost-yang-zuo}, the estimate of $|dh|^2$ should be read as "$|dh|^2\le C$" instead of "$|dh|^2\le C|\log r|^2$".}}.

Geometrically, the finiteness of the energy can be explained as
follows. For any fixed $z_1, \cdots, z_{i-1}, z_{i+1}, \cdots, z_n,
i\le k$, $H_0$ can be considered as a geodesic (and hence harmonic) isometric embedding
\begin{eqnarray*}
&&\exp({\frac 1{2\pi}}\theta_iN_i)\circ\{\exp\big(({\frac 1 2}\log|\log r_i|)Ad_{\exp(\sum_{j\neq i}{\frac 1{2\pi}}\theta_iN_i)}Y_i\big) \\
&&\circ[\exp({\frac 1{2\pi}}\sum_{j\neq i}\theta_iN_i)\circ\exp(\sum_{j\neq i}({\frac 1 2}\log|\log r_i|)Y_i)]\},
\end{eqnarray*}
of a neighborhood of a point at infinity of the upper half plane, say $\{w_i\in \mathbb{C}~|~\Im w_i>\alpha>0\}$ with $w_i=-\sqrt{-1}\log z_i$, into $\mathcal{P}_r$, which is equivariant with respect to $\gamma_i$, as in the $1$-dimensional case \cite{jost-yang-zuo2}. This point is also important in the proof in \S 4; namely, {\bf our local constructions are harmonic on  punctured disks transversal to the divisor}.

Since the metric $ds^2_P$ is a product metric on $(\Delta^*)^k\times\Delta^{n-k}$, for the estimate of the energy density (and energy) of $H_0$, we can consider each $\partial_{z_i}H_0, i\le k$ separately. Again since the map
\[
p: (\{w_i=x_i+\sqrt{-1}y_i\in \mathbb{C}~|~y_i>\alpha>0\}, {\frac{dw_i\wedge d{\overline w}_i}{|{\text Im}w_i|^2}})\to (\Delta^*, {\frac{dz_i\wedge d{\overline z}_i}{|z_i|^2(\log|z_i|)^2}})
\]
where $p(w_i)=z_i=e^{\sqrt{-1}w_i}$, is a Riemannian covering, so restricting to a fundamental domain of $p$, say $\{x_i+\sqrt{-1}y_i\in \mathbb{C}~|~y_i>\alpha>0, 0\le x_i<1\}$, we can write $H_0$ as
\[
H_0=\exp({\frac 1{2\pi}}\sum_{i=1}^k x_iN_i)\circ\exp(\sum_{i=1}^k({\frac 1 2}\log y_i)Y_i);
\]
and estimating $\partial_{z_i}H_0$ is equivalent to estimating
$\partial_{w_i}H_0$. Since for fixed\\
 $w_1, \cdots, w_{i-1}, w_{i+1}, \cdots, w_n$, $H_0$ is a geodesic isometric embedding, we have
\[
\vert\partial_{w_i}H_0\vert={\text{const.}}.
\]
On the other hand, the domain $\{x_i+\sqrt{-1}y_i\in \mathbb{C}~|~y_i>\alpha>0, 0\le x_i<1\}$ has finite volume, so the energy of $H_0$ is finite.
\\
\\
{\bf 3.1.2 Patching local constructions together on a tube neighborhood of the divisor}
Using the above local construction for $H_0$ and a finite partition of
unity on $\overline X$, we can construct a smooth metric on $L_\rho$
which takes the $H_0$ as asymptotic behavior near the divisor $D$,
still denoted by $H_0$. So, using the above complete K\"ahler metric
on $X$, such a metric, as a $\rho$-equivariant map, has finite
energy. However, although our local constructions are harmonic on
punctured disks transversal to the divisor, the metric $H_0$, after
patched together, is not necessarily harmonic on such a punctured
disk. Here, for convenience later on, we do a special patching;
we choose an appropriate partition of unity to patch
these local constructions together so that the metric $H_0$ is
harmonic on some open subsets of the divisor at infinity when
restricted to a certain small punctured disk transversal to the divisor with
the puncture contained in the divisor. This can be done as follows.

For simplicity, we assume in the following that $\dim_\mathbb{C}\overline X=3$, $D=D_1+D_2$ and $D_1\cap D_2\neq\varnothing$; so $\dim_\mathbb{C}D_1\cap D_2=1$. The discussion for the general case is completely similar.

First, we patch together local constructions near the intersection
$D_1\cap D_2$. Take two enough small tube neighborhoods
$U^{12}_{\epsilon'}\subset U^{12}_\epsilon$ of $D_1\cap D_2$ with a
holomorphic projection $\pi^{12}: U^{12}_\epsilon\to D_1\cap D_2$. We
take $U^{12}_\epsilon$ enough small so that for any point $p$ of
$D_1\cap D_2$ there exists a neighborhood $U$ in $D_1\cap D_2$, such
that $(\pi^{12})^{-1}(U)$ can be covered by a local complex coordinate $(z_1, z_2, z_3)$ at $p$ of $\overline X$ with $D_1=\{z_1=0\}$ and $D_2=\{z_2=0\}$. Obviously, the fibres of $\pi^{12}$ are the product of two disks.

Take a {\it finite} open covering $\{U_\alpha\}$ of $D_1\cap D_2$ satisfying that each $U_\alpha$ has the property of the above $U$ and that for each $U_\alpha$ one can take a smaller open subset $U_\alpha'\subset U_\alpha$ so that $U_\alpha'\cap U_\beta'=\varnothing$. Corresponding to $\{U_\alpha'\subset U_\alpha\}$, we can choose a partition of unity $\{\phi_\alpha\}$ with $\phi_\alpha|_{U_\alpha'}\equiv 1$; consequently, $\{\phi_\alpha\circ\pi^{12}\}$ is a partition of unity on $U^{12}_\epsilon$ with $\phi_\alpha\circ\pi^{12}|_{(\pi^{12})^{-1}(U_\alpha')}\equiv 1$. Using this partition of unity, we can patch smoothly all local constructions together along $D_1\cap D_2$ to get a metric of $L_\rho$ on $U^{12}_\epsilon$. We remark that from the above patching process, we can see that the obtained metric is harmonic on $(\pi^{12})^{-1}(U_\alpha')$ when restricted to punctured disks transversal to $D_1$ or $D_2$.

Using the same way, we can patch all local constructions together
along $D_1-U^{12}_{\epsilon'}$ and $D_2-U^{12}_{\epsilon'}$. For
convenience of notations in the next section, we here give some
details. Take two enough small tube neighborhoods
$V^i_{\epsilon'}\subset V^i_\epsilon$ of $D_i-U^{12}_{\epsilon'}$,
$i=1, 2$; for our purposes, w.l.o.g we may assume that $V^i_\epsilon\cap U^{12}_{\epsilon'}$ is empty. Then there is a holomorphic projection $\pi^i: V^i_\epsilon\to D_i-U^{12}_{\epsilon'}$. We also assume that $V^i_\epsilon$ is small enough so that for any interior point $p$ of $D_i-U^{12}_{\epsilon'}$ there exists a neighborhood $V$ in $D_i-U^{12}_{\epsilon'}$, $(\pi^1)^{-1}(V)$ can be covered by a local complex coordinate $(z_1, z_2, z_3)$ at $p$ of $\overline X$ with $D_i=\{z_i=0\}$. Take a {\it finite} open covering $\{V_\alpha\}$ of $D_i-U^{12}_{\epsilon'}$ satisfying that each $V_\alpha$ has the property of the above $V$ and that for each $V_\alpha$ one can take a smaller open subset $V_\alpha'\subset V_\alpha$ so that $V_\alpha'\cap V_\beta'=\varnothing$. Corresponding to $\{V_\alpha'\subset V_\alpha\}$, we can choose a partition of unity $\{\psi_\alpha\}$ with $\psi_\alpha|_{V_\alpha'}\equiv 1$; consequently, $\{\psi_\alpha\circ\pi^{i}\}$ is a partition of unity on $V^{i}_\epsilon$ with $\psi_\alpha\circ\pi^{i}|_{(\pi^{i})^{-1}(V_\alpha')}\equiv 1$. Using this partition of unity, we can patch smoothly all local constructions together along $D_i-U^{12}_{\epsilon'}$ to get a metric of $L_\rho$ on $V^{i}_\epsilon$.

Finally, we can patch smoothly the above three metrics together so
that the metrics on $U^{12}_{\epsilon'},
V^1_\epsilon-U^{12}_{\epsilon}, V^2_\epsilon-U^{12}_{\epsilon}$ are
preserved. Thus, near the divisor, we get a smooth metric which is
harmonic on  certain open subsets, say $(\pi^{12})^{-1}(U_\alpha')\cap
U^{12}_{\epsilon'}$ and $(\pi^{i})^{-1}(V_\alpha')$, when restricted
to a punctured disk transversal to the divisor. We can then extend
the metric smoothly to all of $X$ to get a metric of $L_\rho$, still
denoted by $H_0$.

\subsection{The norm estimate under $H_0$ of a flat section of $L_\rho$}
Here, we want to observe what the asymptotic behavior of the norm of a
flat section of $L_\rho$ under the metric $H_0$ is near the
divisor. To this end, we continue to restrict ourselves to $(\Delta^*)^k\times\Delta^{n-k}$: Since the $\theta_i$-directions have nothing to do with asymptotic behavior of the norm (actually, we can consider $\exp({\frac 1{2\pi}}\sum_{i=1}^k\theta_iN_i)$ as an isometry on $\mathcal{P}_r$), we only need to observe $\exp(\sum_{i=1}^k({\frac 1 2}\log|\log r_i|)Y_i)$. By the previous construction for semisimple elements $Y_i$, they can be diagonalized simultaneously under a suitable basis of $L_\rho$; assuming this and expanding $\exp(\sum_{i=1}^k({\frac 1 2}\log|\log r_i|)Y_i)$, one has the following form
\begin{equation}\label{metric1}
\left(
\begin{array}{ccccc}
\prod_{i=1}^k|\log r_i|^{\frac{a_i^1}{2}}&0&\cdots&0&0\\
0&\prod_{i=1}^k|\log r_i|^{\frac{a_i^2}{2}}&\cdots&0&0\\
\vdots&\vdots&\ddots&\vdots&\vdots\\
0&0&\cdots&\prod_{i=1}^k|\log r_i|^{\frac{a_i^{r-1}
}{2}}&0\\
0&0&\cdots&0&\prod_{i=1}^k|\log r_i|^{\frac{a_i^r}{2}}
\end{array}
\right),
\end{equation}
where
\begin{equation}\label{metric1}
Y_i=\left(
\begin{array}{ccccc}
a_{i}^1&0&\cdots&0&0\\
0&a_{i}^2&\cdots&0&0\\
\vdots&\vdots&\ddots&\vdots&\vdots\\
0&0&\cdots&a_i^{r-1}&0\\
0&0&\cdots&0&a_i^r
\end{array}
\right).
\end{equation}
This gives explicitly the asymptotic behavior of the norm of a flat section
under the metric $H_0$.

\vskip .3cm
\noindent {\bf Remark:} {\it In Cattani-Kaplan-Schmid's theory
\cite{cattani-kaplan-schmid}, the norm estimates of flat sections
under the Hodge metric depend on the order of $N_1, \cdots, N_n$ and use the
notion of weight filtration, and consequently one needs to restrict to the
corresponding sectors of $(\Delta^*)^k$.}

\subsection{The behavior of the differential of $H_0$}
Finally, we also need to understand the asymptotic behavior near the divisor of the differential of $H_0$.
This is very different from the estimates of the norm, where one does not need to consider the $\theta_i$-directions.

First, we do an explicit computation for the differential of $H_0$. We still restrict to $(\Delta^*)^k\times\Delta^{n-k}$,
\[
H_0=\exp({\frac 1{2\pi}}\sum_{i=1}^k\theta_iN_i)\circ\exp(\sum_{i=1}^k({\frac 1 2}\log|\log r_i|)Y_i).
\]
We can consider $\exp({\frac 1{2\pi}}\sum_{i=1}^k\theta_iN_i)$ as an isometry on $\mathcal{P}_r$; so for the ${\bf r}:=(r_1, \cdots, r_n)$-direction, we have
\[
d_{\bf r}H_0=\big(\exp({\frac 1{2\pi}}\sum_{i=1}^k\theta_iN_i)\big)_*\big(\sum_{i=1}^k{\frac{dr_i}{2r_i\log r_i}}Y_i\big).
\]
For the ${\Theta}:=(\theta_1, \cdots, \theta_n)$-direction, if we consider the matter at the identity $I_r$, the differential $d_{\Theta}H_0$ should be read as
\[
{\frac 1{2\pi}}\sum_{i=1}^kN_id\theta_i;
\]
so,
\[
d_{\Theta}H_0=\big(\exp({\frac 1{2\pi}}\sum_{i=1}^k\theta_iN_i)\exp(\sum_{i=1}^k({\frac 1 4}\log|\log r_i|)Y_i)\big)_*{\frac 1{2\pi}}\sum_{i=1}^kN_id\theta_i,
\]
here we consider $H_0$ as $\big(\exp({\frac 1{2\pi}}\sum_{i=1}^k\theta_iN_i)\exp(\sum_{i=1}^k({\frac 1 4}\log|\log r_i|)Y_i)\big)\circ I_r$, denoted it by $\sqrt{H_0}\circ I_r$.
Translating everything into the complex coordinates $(z_1, \cdots, z_n)$, we have
\begin{eqnarray*}
d_{\bf r}H_0&=&{\frac{1}
{4}}\big(\exp({\frac 1{2\pi}}\sum_{i=1}^k\theta_iN_i)\big)_*\sum_{i=1}^k\big({\frac{dz_i}{z_i}}+{\frac{d\overline z_i}{\overline z_i}}\big){\frac{Y_i}{\log|z_i|}},
\\
d_{\Theta}H_0&=&{\frac 1{4\pi\sqrt{-1}}}\big(\sqrt{H_0}\big)_*\sum_{i=1}^kN_i\big({\frac{dz_i}{z_i}}+{\frac{d\overline z_i}{\overline z_i}}\big).
\end{eqnarray*}

As showed before, under the product metric $ds^2_P$ on $(\Delta^*)^k\times\Delta^{n-k}$ and the invariant metric on $\mathcal{P}_r$, $\Vert dH_0\Vert={\text{constant}}$. So, we have
\[
\parallel d_{\bf r}H_0\parallel^2, \parallel d_{\Theta}H_0\parallel^2\le C,
\]
for some positive constant $C$; in particular, $\Vert N_i{\frac{dz_i}{z_i}}\Vert^2\le C.$
Consequently, since $\vert {\frac{dz_i}{z_i}}\vert^2=|\log|z_i||^2$ under the positive metric $ds^2_P$, $N_i$, as an endomorphism of the bundle $L_\rho$, has the following point-wise norm estimate{\footnote{The norm of a tangent vector of $\mathcal{P}_r$ is equivalent to its norm when it is considered as an endomorphism of $L_\rho$, cf. \cite{yang}.}}
\[
\Vert N_i\Vert^2\le C|\log|z_i||^{-2},
\]
where we always consider $\sum_{i=1}^kN_i{\frac{dz_i}{z_i}}$ as a $1$-form homomorphism on $L_\rho$.
\\
\\
{\bf Remark:} {\it We remark that such a norm estimate is not necessarily precise; for example, in the case of $N_i=N_j$, one even has $\Vert N_i\Vert^2\le C|\log|z_i||^{-2}|\log|z_j||^{-2}$.}

\vskip .3cm
We now express the above estimates in terms of weight filtrations as follows. This is easy if we note the relation $[N_i, Y_i]=2N_i$; corresponding to the diagonalisation of $Y_i$, we can consider the weight filtration of $N_i$: $\{W_k\}$ satisfying $N_iW_k\subset W_{k-2}$. If $v$ is a flat section lying in $W_k$, satisfying
$$
\parallel v\parallel_{H_0}^2\sim\prod_{j=1}^k|\log r_j|^{a_j}
$$
and $N_iv\neq 0$, one then has
\[
\parallel N_iv\parallel_{H_0}^2\le C|\log r_i|^{a_i-2}\prod_{j\neq i} |\log r_j|^{a_j}.
\]

By the above computation and argument, we now have
\begin{eqnarray*}
\partial H_0&=&{\frac{1}
{4}}\big(\exp({\frac 1{2\pi}}\sum_{i=1}^k\theta_iN_i)\big)_*\sum_{i=1}^k{\frac{dz_i}{z_i}}{\frac{Y_i}{\log|z_i|}}\\
&&+{\frac 1{4\pi\sqrt{-1}}}\big(\sqrt{H_0}\big)_*\sum_{i=1}^kN_i{\frac{dz_i}{z_i}}
\end{eqnarray*}
and
\begin{prop}
Using ${\frac{dz_i}{z_i}}$ as basis, up to some isometries, $\partial H_0$ is asymptotic to $\sum_{i=1}^kN_i{\frac{dz_i}{z_i}}$ as $|z_i|$ goes to zero. Also, from the previous discussion, we know that when considering $\partial H_0$ as a $1$-form homomorphism on $L_\rho$, it is point-wise bounded under the complete K\"ahler metric on $X$, and hence (in applications) $L^2$-bounded as an operator between certain $L^2$-spaces.
\end{prop}

\section{The proofs of Theorem 1, 2}
In order to deform the initial metric $H_0$ into a harmonic one with the same asymptotic behavior as $H_0$, from now on, we assume that the representation $\rho$ is semisimple. Namely, for any boundary
component $\Sigma$ of $\mathcal{P}_r$, there exists an
element $\gamma\in\pi_1(X)$ satisfying
$\Sigma\cap\rho(\gamma)(\Sigma)=\emptyset$; in other words, the
image of $\rho$ does not fix any boundary component or is not contained in any proper parabolic subgroup \cite{donaldson, corlette, yang}. For similar definitions cf. also \cite{jost-yau, labourie}.

Our strategy is first to deform $H_0$ to a harmonic metric, and then to prove that the harmonic metric has the same asymptotic behavior as $H_0$. The first step is standard with the assumption that the representation $\rho$ is semisimple; we here give a sketch together with some necessary properties.

\subsection{Harmonic metric obtained by deforming the initial metric $H_0$}
In the following, we consider $X=\overline X-D$ as a
complete noncompact manifold with the K\"ahler metric constructed in
Proposition 1; sometimes, we also need the K\"ahler metric to be
locally of the product form $ds^2_P$; we can do this in those cases
that do not involve its derivatives.

Take a sequence of compact
manifolds $\{X_i\}$ (with smooth boundary) such that $X_i\subset
X_{i+1}$ and $\cup_{i=1}^\infty X_i=X$. Then, according to Hamilton
\cite{hamilton} (the theory of Hamilton readily applies to the
equivariant setting), we can find a harmonic metric $H_i$ of $L_\rho$
on $X_i$ with $H_i|_{\partial X_i}=H_0|_{\partial X_i}$ and $E(H_i;
X_i)\le E(H_0; X_i)\le E(H_0)$. Next, we need to prove that there
exists a subsequence of $\{H_i\}$ which converges uniformly on any
compact subset of $X$. Fix a compact subset $X_0$ of $X$ and a point
$p\in X_0$. Then the uniform boundedness of the energies of $\{H_i\}$
imply that the energy densities $e(H_i)$ are also uniformly bounded in
$i$ on $X_0$. Using the semisimplicity of $\rho$ and the uniform
boundedness of the energy densities $e(H_i)$, one can show that
$H_i(p)\in \mathcal{P}_r$ are uniformly bounded in $i$ (here,
considering $H_i$ as equivariant map into $\mathcal{P}_r$,
cf. \cite{yang}); so, $\{H_i(X_0)\subset \mathcal{P}_r\}$ are also
uniformly bounded in $i$. Taking a diagonal sequence, we can find a subsequence of $\{H_i\}$, still denoted by $\{H_i\}$, which converges uniformly on any compact subset of $X$ to a harmonic metric $H$ of $L_\rho$ on $X$. Furthermore, $E(H)\le E(H_0)<\infty$.

Due to the finiteness of energy of the harmonic metric $H$, Siu's Bochner technique for harmonic maps \cite{sampson} (which also applies to the equivariant setting) implies

\begin{lem}
The harmonic metric $H$, as an equivariant map into $\mathcal{P}_r$, is pluriharmonic; namely, when restricted to any complex curve, especially a punctured disk transversal to the divisor $D$, $H$ is harmonic.
\end{lem}

\subsection{$H$ having the same asymptotic behavior as $H_0$}
Now, we can show that $H$ and $H_0$ have the same asymptotic behavior near the divisor; we only need to show that $H$ has the same asymptotic behavior as every local construction $H_0$ in \S3.1. For simplicity of discussion, we continue to restrict to the case that $\dim_\mathbb{C}\overline X=3$, $D=D_1+D_2$ and $D_1\cap D_2\neq\varnothing$ and to use the notations in \S 3.1.2.

We first consider the situation near $D_1\cap D_2$, i.e. in $U^{12}_{\epsilon'}-D_1\cup D_2$. Fix arbitrarily an open set $U_\beta$ in the finite open covering $\{U_\alpha\}$ of $D_1\cap D_2$ and a point $p\in U_\beta'\subset U_\alpha$. By the choice of the open covering, we have a local coordinate $(z_1, z_2, z_3)$ at $p$ of $\overline X$ covering $(\pi^{12})^{-1}(U_\beta')\cap U^{12}_{\epsilon'}$ with $p=(0, 0, 0)$, $D_1=\{z_1=0\}$, and $D_2=\{z_2=0\}$. By the local construction of $H_0$ at $p$ in \S3.1, $H_0$ is harmonic and of finite energy when restricted to a punctured disk (in $(\pi^{12})^{-1}(U_\beta')\cap U^{12}_{\epsilon'}-D_1\cup D_2$) transversal to $D_1$ or $D_2$; on the other hand, $H$ has the same properties by its pluriharmonicity and energy finiteness. Now, we claim

\begin{lem} ${\text{dist}}_{\mathcal{P}_r}(H, H_0)(z_1', z_2', 0)$ {\footnote{Considering $H_, H_0$ as equivariant maps from the universal covering into $\mathcal{P}_r$, the distance function between $H, H_0$ is still equivariant, so can be considered as a function on the base manifold.}} is uniformly bounded for $(z_1', z_2', 0)\in (\pi^{12})^{-1}(U_\beta')\cap U^{12}_{\epsilon'}-D_1\cup D_2$.
\end{lem}
{\it Proof:} First, for a fixed $z_1'$, define
$$
S_\beta^2(z_1')=\{(z_1', z, 0)\in (\pi^{12})^{-1}(U_\beta')\cap U^{12}_{\epsilon'}-D_1\cup D_2\};
$$
w.l.o.g., we can assume it is a punctured disk which is transversal to $D_2$ at $(z_1', 0, 0)$. As pointed out above, $H_0$ and $H$ are harmonic when restricted to $S^2_\beta(z_1')$.
We now prove
\[
\sup_{x\in S^2_\beta(z_1')}{\text{dist}}(H, H_0)(x)\le \sup_{x\in \partial U^{12}_{\epsilon'}\cap S_\beta^2(z_1')}{\text{dist}}(H, H_0)(x).
\]
In order to prove the above inequality, we consider the sequence of harmonic metrics $\tilde H_i$ on $S_\beta^2(z_1')\cap X_i$ with
\[
\tilde H_i|_{\partial U^{12}_{\epsilon'}\cap S_\beta^2(z_1')}=H_i|_{\partial U^{12}_{\epsilon'}\cap S_\beta^2(z_1')}, ~
\tilde H_i|_{\partial X_i\cap S_\beta^2(z_1')}=H_0|_{\partial X_i\cap S_\beta^2(z_1')}.
\]
We remark that $H_i|_{\partial U^{12}_{\epsilon'}\cap S_\beta^2(z_1')}$ converges uniformly to $H|_{\partial U^{12}_{\epsilon'}\cap S_\beta^2(z_1')}$. It is easy to prove that the sequence $\{\tilde H_i\}$ (if necessary, go to a subsequence) converges uniformly on any compact subset of $S_\beta^2(z_1')$ to a harmonic metric $\tilde H$ on $S_\beta^2(z_1')$ and that, by the subharmonicity of ${\text{dist}}(\tilde H_i, H_0)$ {\footnote{ Since $\tilde H_i, H_0$ are harmonic and $\mathcal{P}_r$ has non-positive curvature, a standard computation shows that ${\text{dist}}(\tilde H_i, H_0)$ is subharmonic; cf. e.g. \cite{schoen-yau}.}} on $S_\beta^2(z_1')\cap X_i$ and  the maximum principle,
\[
\sup_{x\in S^2_\beta(z_1')}{\text{dist}}(\tilde H, H_0)(x)\le \sup_{x\in\partial U^{12}_{\epsilon'}\cap S_\beta^2(z_1')}{\text{dist}}(\tilde H, H_0)(x);
\]
and hence $\tilde H$ also has finite energy. Thus, on $S_\beta^2(z_1')$, we have two harmonic metrics $H, \tilde H$ with finite energy and $\tilde H|_{\partial U^{12}_{\epsilon'}\cap S_\beta^2(z_1')}=H|_{\partial U^{12}_{\epsilon'}\cap S_\beta^2(z_1')}$.

We now want to show $\tilde H\equiv H$ on $S_\beta^2(z_1')$, and hence the required inequality.
We observe the distance function ${\text{dist}}(\tilde H, H)$ on $S_\beta^2(z_1')$, denoted by $w$. It is clear that $w$ is subharmonic and $w|_{\partial U^{12}_{\epsilon'}\cap S_\beta^2(z_1')}=0$; on the other hand, by the finiteness of energy of $\tilde H$ and $H$, a simple computation shows that $w$ is of finite energy. Thus, the problem is reduced to show the following lemma (for its proof, cf. Appendix).

\begin{lem} Let $\Delta^*$ be a puncture disk. Assume that $w$ is a non-negative subharmonic function vanishing on the exterior boundary and of finite energy, then $w\equiv 0$.
\end{lem}

Assume that $\sup_{x\in \partial S^2_\beta(z_1')}{\text{dist}}(H, H_0)(x)$ is attained at the point $(z_1', z_2^0, 0)$. Define
\[
S_\beta^1(z_2^0)=\{(z, z_2^0, 0)\in (\pi^{12})^{-1}(U_\beta')\cap U^{12}_{\epsilon'}-D_1\cup D_2\};
\]
also assume that it is a punctured disk which is transversal to $D_1$ at $(0, z_2^0, 0)$.
The same argument as above shows
\[
\sup_{x\in S^1_\beta(z_2^0)}{\text{dist}}(H, H_0)(x)\le \sup_{x\in\partial U^{12}_{\epsilon'}\cap S_\beta^1(z_2^0)}{\text{dist}}(H, H_0)(x).
\]

Assuming that $\sup_{x\in\partial U^{12}_{\epsilon'}\cap S_\beta^1(z_2^0)}{\text{dist}}(H, H_0)(x)$ is attained at the point $(z_1^0, z_2^0, 0)$, then we have
\[
{\text{dist}}(H, H_0)(z_1', z_2', 0)\le {\text{dist}}(H, H_0)(z_1^0, z_2^0, 0)
\]
for $(z_1', z_2', 0)\in (\pi^{12})^{-1}(U_\beta')\cap U^{12}_{\epsilon'}-D_1\cup D_2$. Clearly, we can choose a fixed compact subset of $X$, which is independent of $(z_1', z_2', 0)$ and always contains the corresponding $(z_1^0, z_2^0, 0)$. This proves the lemma.   ~~~~~~~~~~~$\Box$

The arbitrarity of $p$ and the above lemma imply that $H$ and $H_0$ have the same asymptotic behavior at infinity on $U^{12}_{\epsilon'}-D_1\cup D_2$.

Using a similar argument as above, we can also show that $H$ and $H_0$ have the same asymptotic behavior at infinity on $V^{i}_{\epsilon'}-D_i$, i=1, 2. Thus, we have showed that $H$ and $H_0$ have the same asymptotic behavior at infinity on $X$.

We remark that in the above proof, we only use two properties of the harmonic metric $H$: finiteness of energy and pluriharmonicity. So we have
\begin{thm} Any harmonic metric $H$ on $L_\rho$ ($\rho$ is not necessarily semisimple) of finite energy must have the same asymptotic behavior as the local construction {\rm (\ref{H0})} in the introduction near the divisor.
\end{thm}

Now, we turn to the uniqueness problem. This is easy. Assume $H_1, H_2$ are two harmonic metrics with the same asymptotic behavior as $H_0$. Then, ${\text{dist}}_{\mathcal{P}_r}(H_1, H_2)$ is a bounded subharmonic function on $X$. A well-known fact which says that a bounded positive subharmonic function on a complete Riemannian manifold with finite volume has to be constant implies
\[
{\text{dist}}_{\mathcal{P}_r}(H_1, H_2)= {\text{constant}}.
\]
This shows that
$H_1$ and $H_2$
are the same up to an isometry of $\mathcal{P}_r$.

\vskip 1cm
\noindent {\bf Appendix:} Since the punctured disk is conformal to a half cylinder, so Lemma 4 is equivalent to the following

\begin{lem} Let $C$ be a half cylinder with Euclidean metric. Assume that $w$ is a non-negative subharmonic function vanishing on $\partial C$ and of finite energy, then $w\equiv 0$.
\end{lem}

Take a global coordinate $(\theta, y)$ on $C$ with $0\le\theta\le\pi$ and $y\ge 0$; the Euclidean metric is $d\theta^2+dy^2$. Take a sequence of compact subsets $C_n=\{y\le n\}, n=1, 2, \cdots$, and consider a sequence of harmonic functions $u_n$ with $u_n|_{\partial C}=0$ and $u_n|_{y=n}=w_n|_{y=n}$. The maximum principle implies $u\ge w$ on $C_n$ and $E(u_n; C_n)\le E(w; C_n)\le E(w)$. The standard elliptic estimate implies $u_n$ converges to a harmonic function $u$ on $C$ uniformly on any compact subset, which satisfies $u\ge w$, $u|_{\partial C}=0$, and $E(u)\le E(w)$.
So, the proof of Lemma 5 can be reduced to the following fact
which is known to analysts but whose proof we include here for completeness.
\\
\\
{\bf Lemma.} {\it Let $C$ be a half cylinder with Euclidean metric. Assume that $u$ is a non-negative harmonic function with $u|_{\partial C}=0$ and finite energy, then $u\equiv 0.$}
\\
\\
{\it Proof of Lemma.} Take a global coordinate $(\theta, y)$ on $C$ with $0\le\theta\le\pi$ and $y\ge 0$; the Euclidean metric is $d\theta^2+dy^2$. Take a smooth cut-off function $\psi_{r, \rho}$ with $\psi_{r, \rho}(\theta, y)=1, y\le r$, $\psi_{r, \rho}(\theta, y)=0, y\ge r+\rho, \rho> 0$ and $|\nabla\psi_{r, \rho}|\le{\frac 1 \rho}$.

Compute
\begin{eqnarray*}
0&=&-\int_C\Delta u u(\psi_{r, \rho})^2\\
&=&-\int_C\nabla(\nabla u u(\psi_{r, \rho})^2)+\int_C|\nabla u|^2(\psi_{r, \rho})^2+2\int_C\nabla u\psi_{r, \rho}\nabla\psi_{r, \rho}u\\
&\ge&\int_C|\nabla u|^2(\psi_{r, \rho})^2-({\frac 1 2}\int_C|\nabla u|^2(\psi_{r, \rho})^2+2\int_C|\nabla\psi_{r, \rho}|^2u^2).
\end{eqnarray*}
So, we have
\[
\int_C|\nabla u|^2(\psi_{r, \rho})^2\le 4\int_C|\nabla\psi_{r, \rho}|^2u^2.
\]
It is clear that $\int_C|\nabla u|^2(\psi_{r, \rho})^2\to E(u)$, as $r$ goes to infinity. So, if we can choose appropriate $\psi_{r, \rho}$ so that $\int_C|\nabla\psi_{r, \rho}|^2u^2$ can be arbitrarily small, then $E(u)=0$ and hence $u=0$.

Take a sufficiently small $\epsilon$. Since $E(u)<\infty$, $\exists r_0$ so that $\int_{y\ge r}|\nabla u|^2\le\epsilon$ for $r\ge r_0$. On the other hand,
\[
\int_C|\nabla\psi_{r, \rho}|^2u^2\le{\frac 1 {\rho^2}}\int_{r\le y\le r+\rho}u^2.
\]

We now estimate $\int_{r\le y\le r+\rho}u^2$. Since $u$ is harmonic, $|\nabla u|^2$ is subharmonic. So by the average value inequality, we have
\[
|\nabla u|^2(x)\le c\int_{r+i-1\le y\le r+i+2}|\nabla u|^2, ~r+i\le y(x)\le r+i+1,
\]
for $i=0, 1, \cdots, n-1$, where $c$ is a positive constant independent of $r, i$.

First, estimate $|u(\theta, r)-u(\theta, r+n)|$ as follows.
\begin{eqnarray*}
|u(\theta, r)-u(\theta, r+n)|^2&=&|\sum_{i=0}^{n-1}(u(\theta, r+i)-u(\theta, r+i+1))|^2\\
&\le&\Big(\sqrt{c}\sum_{i=0}^{n-1}\sqrt{\int_{r+i-1\le y\le r+i+2}|\nabla u|^2}\Big)^2\\
&\le& 3cn\int_{r-1\le y\le r+n+1}|\nabla u|^2\\
&\le& 3c n\int_{y\ge r_0}|\nabla u|^2\le 3c n\epsilon.
\end{eqnarray*}
In general, we have
\[
|u(\theta, r)-u(\theta, r+s)|^2\le 3c\epsilon s.
\]
On the other hand, we have $|u(\theta, r+s)|^2\le 2|u(\theta, r+s)-u(\theta, r)|^2+2|u(\theta, r)|^2$.

Collecting the above, we have
\begin{eqnarray*}
{\frac 1 {\rho^2}}\int_{r\le y\le r+\rho}u^2&\le& {\frac 1 {\rho^2}}\int_{r\le y\le r+\rho}(2|u(\theta, y)-u(\theta, r)|^2+2|u(\theta, r)|^2)\\
&\le&{\frac 1 {\rho^2}}\int_{r\le y\le r+\rho}(6c\epsilon (y-r)+2|u(\theta, r)|^2).
\end{eqnarray*}
Take $r=r_0+1$ and $\rho$ sufficiently large, we have
\[
{\frac 1 {\rho^2}}\int_{r\le y\le r+\rho}u^2\le c'\epsilon,
\]
for a certain positive constant. Thus, when $ r, \rho$ are sufficiently large, $\int_C|\nabla\psi_{r, \rho}|^2u^2$  becomes arbitrarily small. The proof of the lemma is finished.

\end{document}